\documentclass{birkjour}

\usepackage{amsfonts}
\usepackage{epsfig}
\usepackage{amsmath}
\usepackage{graphicx}
\usepackage{amsmath}
\usepackage{amscd}

 \newtheorem{thm}{Theorem}[section]
 \newtheorem{cor}[thm]{Corollary}
 \newtheorem{lem}[thm]{Lemma}
 \newtheorem{prop}[thm]{Proposition}
 \theoremstyle{definition}
 \newtheorem{defn}[thm]{Definition}
 \newtheorem{nota}[thm]{Notation}
 \theoremstyle{remark}
 
 \newtheorem*{ex}{Example}
 \numberwithin{equation}{section}

\def\bt{\begin{thm}\rm}
\def\et{\end{thm}}
\def\bc{\begin{cor}\rm}
\def\ec{\end{cor}}
\def\bx{\begin{ex}}
\def\ex{\end{ex}}
\def\bl{\begin{lem}\rm}
\def\el{\end{lem}}
\def\bp{\begin{prop}\rm}
\def\ep{\end{prop}}
\def\bn{\begin{defn}\rm}
\def\en{\end{defn}}
\def\bo{\begin{proof}}
\def\eo{\end{proof}}
\def\bj{\begin{nota}\rm}
\def\ej{\end{nota}}

\def\ba{\begin{array}}
\def\ea{\end{array}}
\def\be{\begin{equation}}
\def\ee{\end{equation}}
\def\bd{\begin{description}}
\def\ed{\end{description}}
\def\bu{\begin{enumerate}}
\def\eu{\end{enumerate}}
\def\bi{\begin{itemize}}
\def\ei{\end{itemize}}

\newbox\bigstrutbox
\setbox\bigstrutbox=\hbox{\vrule height12.5pt depth5pt width0pt}
\def\bigstrut{\relax\ifmmode\copy\bigstrutbox\else\unhcopy\bigstrutbox\fi}

\newbox\Bigstrutbox
\setbox\Bigstrutbox=\hbox{\vrule height17.5pt depth5pt width0pt}
\def\Bigstrut{\relax\ifmmode\copy\Bigstrutbox\else\unhcopy\Bigstrutbox\fi}

\def\ds{\displaystyle}

\def\A{{\bf A}}
\def\B{{\bf B}}
\def\C{{\bf C}}
\def\D{{\bf D}}

\def\I{{\bf I}}

\def\L{{\bf L}}
\def\M{{\bf M}}
\def\N{{\bf N}}
\def\O{{\bf O}}

\def\a{{\bf a}}
\def\b{{\bf b}}
\def\c{{\bf c}}
\def\d{{\bf d}}
\def\e{{\bf e}}
\def\f{{\bf f}}

\def\h{{\bf h}}
\def\i{{\bf i}}
\def\j{{\bf j}}
\def\k{{\bf k}}

\def\p{{\bf p}}
\def\q{{\bf q}}

\def\s{{\bf s}}
\def\t{{\bf t}}
\def\u{{\bf u}}
\def\v{{\bf v}}

\def\x{{\bf x}}

\def\lfa{{\acute{\a}}}
\def\lfb{{\acute{\b}}}
\def\lfc{{\acute{\c}}}
\def\lfd{{\acute{\d}}}

\def\lfq{{\acute{\q}}}

\def\lfba{{\acute{\bar{\a}}}}
\def\lfbb{{\acute{\bar{\b}}}}
\def\lfbc{{\acute{\bar{\c}}}}

\def\0{{\bf 0}}
\def\1{{\bf 1}}
\def\2{{\bf 2}}
\def\3{{\bf 3}}
\def\4{{\bf 4}}
\def\5{{\bf 5}}
\def\6{{\bf 6}}
\def\7{{\bf 7}}
\def\8{{\bf 8}}
\def\9{{\bf 9}}

\def\bq{\overline{\q}}

\def\qh{\mathbb{H}}

\def\overlinea{\overline{\a}}

\begin{document}

\title[Basis-free Solution to General Linear Quaternionic Equation]
{Basis-free Solution to General Linear Quaternionic Equation}

\author[Changpeng Shao, Hongbo Li, Lei Huang]
{Changpeng Shao, Hongbo Li, Lei Huang}

\address{
KLMM, AMSS and UCAS\\
Chinese Academy of Sciences\\
Beijing 100190, China}

\email{\\
shaochangpeng11@mails.ucas.ac.cn,\\ hli@mmrc.iss.ac.cn,\\ lhuang@mmrc.iss.ac.cn}



%\date{January 1, 2004}

\begin{abstract}
A linear quaternionic equation in one quaternionic variable $\q$ is of the form 
$\a_1\q\b_1+\a_2\q\b_2+\cdots +\a_m\q\b_m=\c$,
where the $\a_i, \b_j, \c$ are given quaternionic coefficients. If introducing basis 
elements $\i, \j, \k$ of pure quaternions,
then the quaternionic equation becomes four linear equations in four unknowns over the 
reals, and solving such equations is trivial.
On the other hand, finding a quaternionic rational function expression of the 
solution that involves
only the input quaternionic coefficients and their conjugates, called a 
{\it basis-free solution}, is non-trivial.

In 1884, Sylvester initiated the study of basis-free solution to linear quaternionic equation. 
He considered the three-termed
equation $\a\q+\q\b=\c$, and found its solution 
$\q=(\a^2+\b\overline{\b}+\a(\b+\overline{\b}))^{-1}(\a\c+\c\overline{\b})$
by successive left and right multiplications.
In 2013, Schwartz extended the technique to the four-termed equation, and obtained the basis-free 
solution in explicit form.

This paper solves the general problem for arbitrary number of terms in the non-degenerate 
case. 
\end{abstract}

\maketitle

\section{Introduction}
\setcounter{equation}{0}

A {\it quaternionic variable} is $\q:=u+x\i+y\j+z\k$,
where the $u,x,y,z$ are real-valued variables,
and $1, \i,\j,\k$ are the fixed basis of quaternions. 
A {\it quaternionic monomial} of degree
$m$ is of the form $\a_1\q\a_2\q\cdots\a_{m}\q\a_{m+1}$, where the quaternionic 
variable $\q$ occurs $m$ times, and the $\a_l$ are quaternionic coefficients.

By rewriting the quaternionic 
variable and the quaternionic coefficients as
linear combinations of the fixed basis, a quaternionic polynomial
equation is converted into four real polynomial equations in four real variables, 
and methods of real polynomial system
solving can be applied to find solutions for the $u,x,y,z$, then 
for the unknown $\q$.
So if the basis elements $\i, \j, \k$ are introduced, then solving quaternionic 
equations is reduced to solving real polynomial equations.

For example, for a
{\it linear quaternionic equations} of $n$ terms \cite{rodman}, {\it i.e.}, an equation of the form
\be
\c_1\q\b_1+\c_2\q\b_2+\cdots+\c_{n-1}\q\b_{n-1}=\d,
\label{eqn:general}
\ee
by introducing real-valued variables $u,x,y,z$ for quaternionic variable
$\q$,
and real-valued coordinates for the $\c_i, \b_j, \d$, (\ref{eqn:general}) is converted into 
four real linear equations in four unknowns, called the {\it associated real linear system}
of the equation. Standard methods such as Gaussian elimination can be used to solve the linear system.

For an input quaternionic equation where the basis elements $\i, \j, \k$ do not occur explicitly
in the coefficients, to keep the equation solving within the framework of quaternions instead of 
real numbers, it is sometimes desired that the quaternionic solution be expressed as
a quaternionic rational function in the input quaternionic coefficients and their conjugates only, 
so that if the input equation does not involve the basis elements $\i,\j,\k$, nor does the
solution. Such a solution, if exists, is 
called a {\it basis-free solution}.

For example, for the linear quaternionic equation (\ref{eqn:general}), if $n=2$, then
solving $\c\q\b=\d$ is trivial, and the basis-free solution is 
$\q=\c^{-1}\d\b^{-1}$ under the non-degeneracy condition
$\c\b\neq 0$. If $n=3$, by multiplying both sides of (\ref{eqn:general}) from the right 
with $\b_1^{-1}$
and from the left with $\c_2^{-1}$, we get the following {\it Sylvester's equation} 
\cite{sylvester}:
\be
\label{sylvester-equation}
\s\q+\q\t=\u,
\ee
where $\s=\c_2^{-1}\c_1$, \ $\t=\b_2\b_1^{-1}$, and $\u=\c_2^{-1}\d\b_1^{-1}$ under 
the non-degeneracy
condition $\c_2\b_1\neq 0$. It is well known that the basis-free solution of (\ref{sylvester-equation}) 
can be obtained as follows:
\bu
\item Multiply both sides of (\ref{sylvester-equation}) from the left with $\s$:
\be \label{step1-two-terms}
\s^2\q+\s\q\t=\s\u.
\ee

\item Multiply both sides of (\ref{sylvester-equation}) from the right with $\overline{\t}$:
\be \label{step2-two-terms}
\s\q\overline{\t}+\t\overline{\t}\q=\u\overline{\t}
\ee

\item Add up (\ref{step1-two-terms}) and (\ref{step2-two-terms}):
\be
(\s^2+\t\overline{\t}+(\t+\overline{\t})\s)\q=\s\u+\u\overline{\t}.
\ee

\item
The solution is
\be
\q=(\s^2+\t\overline{\t}+(\t+\overline{\t})\s)^{-1}(\s\u+\u\overline{\t})
\label{soln:syl}
\ee
under the non-degeneracy condition $\s^2+\t\overline{\t}+(\t+\overline{\t})\s\neq 0$.
\eu

{\bf Notation.}\
Following \cite{deleo}, for arbitrary quaternions $\u, \v$,
we use ``$(\u\,|\,\v)$" to denote the following linear operator upon quaternions:
\be
(\u\,|\,\v)\q:=\u\q\v, \hbox{ for any }\q\in \qh.
\ee
This operator satisfies the following composition rule:
\be
(\u_1\,|\,\v_1)(\u_2\,|\,\v_2)=(\u_1\u_2\,|\,\v_2\v_1).
\ee

By the above notation, the procedure of solving (\ref{sylvester-equation}) can be 
written succinctly as the following identity:
\be
\Big(
(\s\,|\,1)+(1\,|\,\overline{\t})\Big)
\Big(
(\s\,|\,1)+(1\,|\,\t)\Big)
=\big(\s^2+\t\overline{\t}+(\t+\overline{\t})\s\,|\,1\big).
\label{proc:sylvester}
\ee

In the literature, 
\cite{johnson} extended Sylvester's equation to algebraic division ring, and
\cite{helmstetter, porter} classified the solutions of Sylvester's 
equation in various degenerate cases. 
For the more general linear quaternionic equation
(\ref{eqn:general}), \cite{janovska, janovska1, shpakivskyi} investigated the 
non-degeneracy conditions of the equation by considering its associated real linear system.
When the number of terms $n=4$, \cite{schwartz} obtained the basis-free solution
in explicit form in the non-degenerate case, by reducing the number of terms to 
three with successive left and right multiplications.

There are also many results on properties of nonlinear quaternionic equations and
algorithms to solve them
\cite{eilenberg, janovska0, janovska2, kalantari, niven41, pogorui, serodio, turner}. These results, 
concerning only the derived real polynomial system of the input equations,
do not lead to basis-free solution. Even in the simplest case of linear quaternionic equation
(\ref{eqn:general}) whose associated real linear system is non-degenerate, finding the 
basis-free solution is highly non-trivial, and 
there is still no result for $n>4$ in the literature.

In \cite{schwartz}, the following idea is proposed to obtain a 
non basis-free solution of
(\ref{eqn:general}). Introduce the coordinate form of the quaternionic coefficients 
$\b_p$ for $p=1..n-1$:
\be
\b_p=b_{p0}+b_{p1}\i+b_{p2}\j+b_{p3}\k.
\ee
Substituting them into (\ref{eqn:general}) and expanding the result, we get
\be
\a_0\q+ \a_1\q\i + \a_2\q\j + \a_3\q\k = \d,
\label{rev:start}
\ee
where for $i=0..3$, 
\be
\a_i=\sum_{p=1}^{n-1} b_{pi}\c_p\in \qh. 
\label{expr:al}
\ee
(\ref{rev:start}) 
is called the {\it revised starting form} \cite{schwartz}
of equation (\ref{eqn:general}).
Non basis-free solution in explicit form can be computed by successive left and 
right multiplications, or simply by Gaussian elimination.

To find the basis-free solution of (\ref{eqn:general}), our first approach is first to 
derive the explicit form of the non basis-free solution of (\ref{eqn:general}) by
solving (\ref{rev:start}), then to convert it into a basis-free quaternionic expression. 
This procedure is difficult, but we manage to make it by brute-force
try and error. The outline of this approach is as follows.

Let $\f$ be the $\mathbb R$-linear isomorphism from $\qh$ to
${\mathbb R}^4$ mapping $1,\i,\j,\k$ to the orthonormal basis $\e_0, \e_1, \e_2, \e_3$
of ${\mathbb R}^4$. The associated real linear system of (\ref{rev:start}) is of the form
$\A\f(\q)=\f(\d)$, where $\A$ is a $4\times 4$ real matrix whose entries are polynomials
in the coordinates of the $\a_l$ for $l=0..3$. 
Under the non-degeneracy condition $\det(\A)\neq 0$, 
the solution to (\ref{rev:start}) is $\f(\q)=\A^{-1}\f(\d)
=adj(\A) \f(\d)/\det(\A)$, where $adj(\A)$ is the adjugate matrix of $\A$. 
The result is a basis-dependent vector-valued rational function in the 
coordinates of the $\a_l$ and $\d$, where each numerator is a 1984-termed polynomial, 
and the denominator $\det(\A)$ is a 196-termed polynomial.

By setting up a list of correspondences between several Clifford algebraic expressions
over ${\mathbb R}^4$ and their basis-free preimages in $\qh$ under $\f$, we manage to
rewrite $adj(\A)$ as the matrix form of the following linear operator acting upon $\qh$:
$(\p_0\,|\,1)+(\p_1\,|\,\i)+(\p_2\,|\,\j)+(\p_3\,|\,\k)$,
where each $\p_l$ is a basis-free quaternionic polynomial in the $\a_l$ and their conjugates.
Similarly, we rewrite $\det(\A)$ as a basis-free quaternionic polynomial in the $\a_l$ 
and their conjugates. 

Then by
substituting the expressions (\ref{expr:al}) of the $\a_l$ into the above result, after 
some hard manipulations upon the expanded summations, we get the basis-free 
quaternionic solution to (\ref{eqn:general}).
As a byproduct, we get a procedure of left and
right multiplications by expressions of the coefficients of (\ref{eqn:general}), so that the
basis-free solution is obtained by adding up the results of the multiplications, 
similar to the procedure from (\ref{step1-two-terms}) to (\ref{soln:syl})
in solving Sylvester's equation.

The first approach is not only complicated, but also unable to 
interpret why some Clifford algebraic expressions over ${\mathbb R}^4$ occur naturally in
the solution, and why they can significantly reduce the size of the solution expression. 

Our second approach is first to embed the algebra of quaternions into 
$CL(\mathbb{R}^4)$, so that the input equation becomes a linear Clifford 
algebraic equation, then to solve the linear Clifford 
algebraic equation, and then to project the solution onto $\qh$ to get  
the basis-free solution of the input quaternionic equation. It turns out that this approach
is much easier, and very elegant. 

There are various embeddings
of $\qh$ as a subalgebra into $CL(\mathbb{R}^4)$. A
classical algebraic homomorphism \cite{hestenes} is generated by the correspondence 
between the basis elements 
$1\leftrightarrow 1,\, -\e_{23}\leftrightarrow \i,\, \e_{13}\leftrightarrow \j$, and $-\e_{12}\leftrightarrow \k$.
To include both this algebraic homomorphism and the linear isomorphism 
$\f: \qh\longrightarrow {\mathbb R}^4$ within the same framework, 
we need a linear map $\pi$ from $CL(\mathbb{R}^4)$ to $\qh$, so that $\pi$ restricted to
$\mathbb{R}^4$ is just $\f^{-1}$, while $\pi$ restricted to the subalgebra 
$\langle 1, \e_{12}, \e_{13}, \e_{23} \rangle_{\mathbb R}$ gives the classical 
algebraic homomorphism. Such a linear map can be easily constructed, and it 
``lift"s the input
quaternionic equation (\ref{eqn:general}) to the following linear 
Clifford algebraic equation:
\be
\sum_{r=1}^{n-1}\f(\c_r)\f(\bar{\q})\f(\b_r)(1+\e_{0123})=\f(\d)(1+\e_{0123}).
\ee
Solving this equation can be done easily by hand. The solution when projected to 
$\qh$ by $\pi$ gives the solution that is obtained by the first approach in
exactly the same form.

Moreover,
once equation (\ref{eqn:general}) is solved, then a more general linear quaternionic equation involving 
one quaternionic variable and its conjugate can be solved readily to 
obtain a basis-free solution.

So this paper solves the long-lasting problem of solving the linear quaternionic 
equation of arbitrary number of terms for the basis-free solution in explicit form, 
by a procedure of successive left and right multiplications with the input
quaternionic coefficients 
and their conjugates, in the non-degenerate case.
The content is arranged as follows. 
Section 2 introduces some notations and basics of linear quaternionic equations. 
Section 3 establishes the correspondence between 
some typical Clifford algebraic expressions over 
$\mathbb{R}^4$ and their basis-free quaternionic representations.
Section 4 presents the associated real linear system approach.
Section 5 presents the Clifford algebra approach. Section 6 solves
the general linear quaternionic equation with conjugate in the non-degenerate case.

\section{Quaternions and linear quaternionic equations}
\setcounter{equation}{0}

In a quaternionic variable $\q:=u+x\i+y\j+z\k$, $u$ is called the {\it scalar part} (or real part), 
and
$x\i+y\j+z\k$ is called the {\it vector part} (or pure imaginary part). The 
quaternionic conjugate is denoted by the
over-bar symbol, while the $\mathbb R$-linear operators extracting the real 
part and the pure imaginary part are denoted by ``$\textmd{Re}$" and 
``$\textmd{Im}$" respectively:
\be
\bq=u-x\i-y\j-z\k, \hskip .4cm
\textmd{Re}(\q)=\frac{\q+\bq}{2},\hskip .4cm
\textmd{Im}(\q)=\frac{\q-\bq}{2}.
\ee
Obviously
$\q\bq=\bq\q=u^2+x^2+y^2+z^2$, and $\q+\bq=2u$; they
are both real-valued.

%Although variable $u$ can be expressed by $\q$ and $\overline{\q}$,
%the other three variables $x,y,z$ cannot be expressed without the introduction of 
%the basis elements $\i,\j,\k$.
%With the three basis elements, the four real variables have the following expressions:
%\be\ba{lll}
%4u      &=& \ds \q-\i\q\i-\j\q\j-\k\q\k;  \\ [1mm]
%4x      &=& \ds \j\q\k-\k\q\j-\i\q-\q\i; \\ [1mm]
%4y      &=& \ds \k\q\i-\i\q\k-\j\q-\q\j; \\ [1mm]
%4z      &=& \ds \i\q\j-\j\q\i-\k\q-\q\k.
%\ea\ee

Let $\e_0, \e_1, \e_2, \e_3$ be a fixed basis of $\mathbb{R}^4$.
Define a linear isomorphism from $\qh$ to $\mathbb{R}^4$ as following:
\be
\f: \hskip .4cm
1\mapsto \e_0, \hskip .2cm
\i\mapsto \e_1, \hskip .2cm
\j\mapsto \e_2, \hskip .2cm
\k\mapsto \e_3.
\ee
For any $\a\in \qh$, denote 
\be
\lfa:=\f(\a),\ \ \
\lfba:=\f(\bar{\a}).
\ee

%{\it Remark}: Alternatively we can define linear isomorphisms from $\qh$ to $\mathbb{R}^{3,1}$
%or $\mathbb{R}^{1,3}$; they lead to the same result for the solving of the linear quaternionic
%equation. For simplicity, in this paper we only
%use $\mathbb{R}^4$ and the Clifford algebra $CL(\mathbb{R}^4)$ it generates.

Under the above linear isomorphism between $\qh$ and ${\mathbb R}^4$, 
the left (or right) multiplication by $\p=
p_0+p_1\i+p_2\j+p_3\k\in \qh$ can be taken as a linear transformation in ${\mathbb R}^4$.
The two multiplications have the following $4\times 4$ real matrix form:
\be
(\p\,|\,1)=\left(
          \begin{array}{rrrr}
            p_0 & -p_1 & -p_2  & -p_3 \\
            p_1 &  p_0 & -p_3  & p_2 \\
            p_2 &  p_3 & p_0   & -p_1 \\
            p_3 & -p_2 & p_1   & p_0 \\
          \end{array}
        \right),
        \hskip .4cm
(1\,|\,\p)=\left(
          \begin{array}{rrrr}
            p_0 & -p_1       & -p_2       & -p_3 \\
            p_1 & \hfill p_0 & \hfill p_3 & -p_2 \\
            p_2 & -p_3       & \hfill p_0 & \hfill p_1 \\
            p_3 & \hfill p_2 & -p_1       & \hfill p_0 \\
          \end{array}
        \right).
\label{a:matrix}
\ee
They have the following simple properties:
\be\ba{ll}
\hbox{(i).} & \hskip 1cm
(\overline{\p}\,|\,1)=(\p\,|\,1)^T;\\

\hbox{(ii).} & \hskip 1cm
(1\,|\,\overline{\p})=(1\,|\,\p)^T;\\

\hbox{(iii).} & \hskip 1cm
\det(\p\,|\,1)=\det(1\,|\,\p)=(\p\overline{\p})^2.

\ea
\ee

For any unit quaternion $\p$, both $(\p\,|\,1)$ and $(1\,|\,\p)$ are in $SO(4)$. Conversely,
any element $\M\in SO(4)$ is of the form
$\M=(\p\,|\,1)(1\,|\,\q)$ for two unit quaternions $\p,\q$ \cite{elfrinkhof}.

Consider the revised starting form (\ref{rev:start}). In matrix form,
\be\ba{lrl}
\A &:=&
(\a_0\,|\,1)+(\a_1\,|\,\i)+(\a_2\,|\,\j)+(\a_3\,|\,\k)\\
[1mm]
&=&
\left(\begin{array}{rr}
       a_{00}-a_{11}-a_{22}-a_{33} & -a_{01}-a_{10}-a_{23}+a_{32}   \\
       a_{01}+a_{10}-a_{23}+a_{32} & a_{00}-a_{11}+a_{22}+a_{33}\\
       a_{02}+a_{13}+a_{20}-a_{31} & a_{03}-a_{12}-a_{21}-a_{30}\\
       a_{03}-a_{12}+a_{21}+a_{30} & -a_{02}-a_{13}+a_{20}-a_{31}
\ea
\right. \\

&& \hskip 2cm
\left.\ba{rr}
        -a_{02}+a_{13}-a_{20}-a_{31} & -a_{03}-a_{12}+a_{21}-a_{30} \bigstrut\\
       -a_{03}-a_{12}-a_{21}+a_{30}  &  a_{02}-a_{13}-a_{20}-a_{31}  \\
         a_{00}+a_{11}-a_{22}+a_{33} & -a_{01}+a_{10}-a_{23}-a_{32}  \\
         a_{01}-a_{10}-a_{23}-a_{32} &  a_{00}+a_{11}+a_{22}-a_{33}
\end{array}\right).
\label{express:A}
\ea
\ee

\bp
For any $4\times 4$ real matrix $\M$, there exist $\p_l\in \qh$ for $l=0..3$ such that
\be
\M=(\p_0\,|\,1)+(\p_1\,|\,\i)+(\p_2\,|\,\j)+(\p_3\,|\,\k).
\ee
More explicitly, if $\M=(m_{ij})_{i,j=0..3}$, then
\be\ba{lll}
4\p_0 &=& \phantom{-} (m_{00}+m_{11}+m_{22}+m_{33})-(m_{01}-m_{10}+m_{23}-m_{32})\i \\
      & & -(m_{02}-m_{13}-m_{20}+m_{31})\j-(m_{03}+m_{12}-m_{21}-m_{30})\k, \\ [1.5mm]
4\p_1 &=& -(m_{01}-m_{10}-m_{23}+m_{32})-(m_{00}+m_{11}-m_{22}-m_{33})\i \\
      & & -(m_{03}+m_{12}+m_{21}+m_{30})\j+(m_{02}-m_{13}+m_{20}-m_{31})\k, \\ [1.5mm]
4\p_2 &=& -(m_{02}+m_{13}-m_{20}-m_{13})+(m_{03}-m_{12}-m_{21}+m_{30})\i \\
      & & -(m_{00}-m_{11}+m_{22}-m_{33})\j-(m_{01}+m_{10}+m_{23}+m_{32})\k, \\ [1.5mm]
4\p_3 &=& -(m_{03}-m_{12}+m_{21}-m_{30})-(m_{02}+m_{13}+m_{20}+m_{31})\i \\
      & & +(m_{01}+m_{10}-m_{23}-m_{32})\j-(m_{00}-m_{11}-m_{22}+m_{33})\k.
\ea
\label{soln:a4}
\ee
As a corollary, any 4 real linear equations in 4 unknowns can be written as a linear quaternionic equation.
\label{matrix-quaternion-form}
\ep

\bo We need to solve the equation $\A=\M$ for the unknowns $a_{ij}$ from the expression of $\A$
in (\ref{express:A}). The entries of $\A$ can be classified into four groups: (i) the four diagonal elements,
as they are linear expressions in $a_{00}, a_{11}, a_{22}, a_{33}$;
(ii) the four anti-diagonal elements (those with row and column numbers $(1,4), (2,3), (3,2), (4,1)$), as they are
linear expressions in $a_{03}, a_{30}, a_{12}, a_{21}$; (iii) the four elements
with row and column numbers $(1,2)$, $(2,1)$, $(3,4)$, $(4,3)$, as they are
linear expressions in $a_{01}, a_{10}, a_{23}, a_{32}$; (iv) the four elements
with row and column numbers $(1,3), (3,1), (2,4), (4,2)$, as they are
linear expressions in $a_{02}, a_{20}, a_{13}, a_{31}$.

First consider group (i). The four equalities of the four corresponding elements in $\A$ and $\M$ respectively can be written as
\be\left(
  \begin{array}{rrrr}
    1 & -1       & -1       & -1 \\
    1 & -1       & \hfill 1 & \hfill 1 \\
    1 & \hfill 1 & -1       & \hfill 1 \\
    1 & \hfill 1 & \hfill 1 & -1 \\
  \end{array}
\right)\left(
         \begin{array}{c}
           a_{00} \\
           a_{11} \\
           a_{22} \\
           a_{33} \\
         \end{array}
       \right)
=\left(\ba{c}
m_{00} \\
m_{11} \\
m_{22} \\
m_{33}
\ea
\right).
\ee
The determinant of the coefficient matrix is nonzero, so the equations have a unique solution for the unknowns
$a_{00}, a_{11}, a_{22}, a_{33}$.

Similarly, in each of the other three groups,
the four equalities of the four corresponding elements in $\A$ and $\M$
respectively form a non-degenerate linear system, so the solutions
are unique. Solving all the linear equations, we get the solution (\ref{soln:a4}).
\eo

\section{Quaternionic representations of some Clifford algebraic expressions over ${\mathbb R}^4$}
\setcounter{equation}{0}

The {\it Clifford algebra} over $\mathbb{R}^4$ is the quotient of the tensor algebra
generated by $\mathbb{R}^4$ modulo the two-sided ideal generated by elements of the form
$\x\otimes \x-\x\cdot \x$ for all vectors $\x\in \mathbb{R}^4$, where the dot symbol denotes the
inner product in $\mathbb{R}^4$. The associative product induced by the quotient from
the tensor product is called the {\it Clifford product}, denoted by the juxtaposition of elements.
The antisymmetrization of the Clifford product is called the {\it exterior product}, also 
known as the
{\it outer product}, and is denoted by the wedge symbol. 

Clifford algebra $CL(\mathbb{R}^4)$ is ${\mathbb Z}_2$-graded. The Clifford product of even number of
vectors is called an {\it even element}, so is any linear combination of even elements. Similarly, the
Clifford product of odd number of
vectors is called an {\it odd element}, so is any linear combination of odd elements.
The Clifford algebra is also ${\mathbb Z}_5$-graded, with grades ranging from 0 to 4: a 0-graded
element is a scalar, a 1-graded element is a vector, and an $r$-graded element for $r>1$ is the 
linear
combination of outer products of $r$ vectors.
The outer product of vectors in $\mathbb{R}^4$ whose number is $>4$
is zero. A general element $\A$ of $CL(\mathbb{R}^4)$ is a linear combination of its
{\it $i$-graded parts},
for $i$ ranging from 0 to 4; the {\it $i$-graded part} of $\A$ is denoted by $\langle \A \rangle_i$.

\bl \cite{li}
For any vectors $\x_i\in \mathbb{R}^4$,
\be\ba{lll}
2 \langle \x_1\x_2\cdots  \x_{2k+1}\rangle_1 &=&
\x_1\x_2\cdots  \x_{2k+1}+\x_{2k+1}\cdots\x_2\x_1, \\

2 \langle \x_1\x_2\cdots  \x_{2k+1}\rangle_3 &=&
\x_1\x_2\cdots  \x_{2k+1}-\x_{2k+1}\cdots\x_2\x_1, \\

2 \langle \x_1\x_2\cdots \x_{2k}\rangle_2 &=&
\x_1\x_2\cdots  \x_{2k}-\x_{2k}\cdots\x_2\x_1, \\

4 \langle \x_1\x_2\cdots \x_{2k}\rangle_4 &=&
\x_1\x_2\cdots  \x_{2k}+\x_{2k}\cdots\x_2\x_1 \\

&& \hfill
-\x_{2k}\x_1\x_2\cdots  \x_{2k-1}
-\x_{2k-1}\cdots\x_2\x_1\x_{2k}, \\

4 \langle \x_1\x_2\cdots \x_{2k}\rangle_0 &=&
\x_1\x_2\cdots  \x_{2k}+\x_{2k}\cdots\x_2\x_1 \\

&& \hfill
+\x_{2k}\x_1\x_2\cdots  \x_{2k-1}
+\x_{2k-1}\cdots\x_2\x_1\x_{2k}.
\ea
\ee
\label{lem:ungrade}
\el

The {\it reverse} ``$^\dagger$" in $CL(\mathbb{R}^4)$ is a linear operator defined as follows:
for any $\x_1\x_2\cdots \x_r$ where $\x_i\in \mathbb{R}^4$,
$(\x_1\x_2\cdots \x_r)^\dagger:=\x_r\cdots \x_2\x_1$.
The {\it inverse} of an element $\A\in CL(\mathbb{R}^4)$, denoted by $\A^{-1}$, is an element
of $CL(\mathbb{R}^4)$ satisfying $\A\A^{-1}=1$. Not every element of $CL(\mathbb{R}^4)$ has inverse.
For example, a vector of $\mathbb{R}^4$ is invertible if and only if it is nonzero.

%The graded adjoint action $Ad^*_{\y}$ of $\y$ upon $CL(\mathbb{R}^4)$ is a linear operator defined as follows:
%for any $\x_1, \ldots \x_r\in \mathbb{R}^4$,
%\be
%Ad^*_{\y}(\x_1\x_2 \cdots \x_r)=(-1)^r \y^{-1}(\x_1\x_2 \cdots \x_r)\y.
%\ee
%Obviously $Ad^*_{\y}(\x_1\x_2 \cdots \x_r)=Ad^*_{\y}(\x_1)\, Ad^*_{\y}(\x_2) \cdots Ad^*_{\y}(\x_r)$.

For fixed orthonormal basis $\e_0, \e_1, \e_2, \e_3$ of $\mathbb{R}^4$, denote 
$\e_{pqr}:=\e_p\e_q\e_r$.
Obviously $\e_{pqr}=\e_p\wedge \e_q\wedge \e_r$. The induced basis of $CL(\mathbb{R}^4)$ is
\be
1, \e_0, \e_1, \e_2, \e_3, \e_{01}, \e_{02}, \e_{03}, \e_{23}, \e_{13}, \e_{12},
\e_{012}, \e_{013}, \e_{023}, \e_{123}, \e_{0123}.
\ee
We usually denote $\I_4:=\e_{0123}$. It is easy to see that $\I_4^\dagger=\I_4^{-1}=\I_4$,
and for any $\x\in \mathbb{R}^4$,
\be
\x\I_4=-\I_4\x.
\ee

In $CL(\mathbb{R}^4)$, the {\it conjugate} is a linear operator defined as follows:
for any $\A\in CL(\mathbb{R}^4)$, its conjugate is
\be
\overline{\A}:=\e_0\A\e_0.
\ee
For example, $\overline{\I_4}=-\I_4$, and for any $\x=\sum_{i=0}^3 x_i \e_i\in \mathbb{R}^4$,
$\overline{\x}=x_0\e_0-\sum_{i=1}^3 x_i \e_i$. Obviously $\overline{\A\B}=
\overline{\A}\,\overline{\B}$ for any $\A,\B\in CL(\mathbb{R}^4)$,
and $\lfba=\overline{\lfa}$ for any $\a\in \qh$. 

The {\it dual} of an element $\A\in CL(\mathbb{R}^4)$, denoted by $\A^\sim$, is defined by
$\A^\sim:=\A\I_4$. It is independent of the choice of the orthonormal basis
$\e_0, \e_1, \e_2$, $\e_3$ of $\mathbb{R}^4$. The dual of $\x_1\wedge \x_2\wedge \x_3\wedge \x_4$ for
$\x_i\in \mathbb{R}^4$ is usually denoted by $[\x_1\x_2\x_3\x_4]$, called the {\it bracket} of
the four vectors.

For example, given $\a_l=u_l+x_l\i+y_l\j+z_l\k\in \qh$ for $l=1..4$, then
\be
[\lfa_1 \lfa_2 \lfa_3 \lfa_4]=
\left|\ba{llll}
u_1 & x_1 & y_1 & z_1 \\
u_2 & x_2 & y_2 & z_2 \\
u_3 & x_3 & y_3 & z_3 \\
u_4 & x_4 & y_4 & z_4
\ea\right|.
\ee
By Lemma \ref{lem:ungrade},
the bracket has a 4-termed representation in Clifford algebraic polynomial form. 

Now we establish the shortest quaternionic polynomial forms of some
typical Clifford algebraic expressions in $CL(\mathbb{R}^4)$, whose
vector variables are images of quaternionic variables under $\f$. 

\bp
For any $\a_1, \a_2\in \qh$, 
\be
\lfa_1\cdot \lfa_2=\lfba_1\cdot  \lfba_2
=\frac{\a_1 \overline{\a_2}+\a_2 \overline{\a_1}}{2}.
\ee
\ep

\bo
For any $\a_l=u_l+x_l\i+y_l\j+z_l\k\in \qh$ where $l=1,2$,
$
\lfa_1\cdot  \lfa_2=\lfba_1\cdot  \lfba_2
=u_1u_2+x_1x_2+y_1y_2+z_1z_2
$
is the inner product in $\mathbb{R}^4$. The conclusion follows this and
$
\textmd{Re}(\a_1\overline{\a_2}) = \textmd{Re}(\a_1)\, \textmd{Re}(\a_2)+\f(\textmd{Im}(\a_1))\cdot \f(\textmd{Im}(\a_2)).
$
\eo

\bl
For any $\a_l\in \qh$ where $l=1..4$,
\be
[\lfa_1 \lfa_2 \lfa_3 \lfa_4]
=- [\lfba_1 \lfba_2 \lfba_3 \lfba_4].
\ee
\el

\bo
The conjugate operator defines an orthogonal transformation
of determinant $-1$ in ${\mathbb R}^4$, so
\[
-\lfa_1 \wedge \lfa_2 \wedge  \lfa_3 \wedge  \lfa_4=
\overline{\lfa_1 \wedge \lfa_2 \wedge  \lfa_3 \wedge  \lfa_4}
=\overline{\lfa_1} \wedge \overline{\lfa_2} \wedge  \overline{\lfa_3} \wedge  \overline{\lfa_4}
=\lfba_1 \wedge \lfba_2 \wedge \lfba_3 \wedge \lfba_4.
\]

\eo

For any $\a_l\in \qh$, define the following elements in $\qh$:
\be
\label{note:check}
\ba{lll}
\a_{pqr}^\sim &:=& \f^{-1}\big( (\lfa_p\wedge \lfa_q\wedge \lfa_r)^\sim\big), \\

\overline{\a}_{pqr}^\sim &:=& \f^{-1}\big((\lfba_p\wedge
\lfba_q\wedge \lfba_r)^\sim\big).\bigstrut
\ea
\ee
By Lemma \ref{lem:ungrade}, 
$\lfa_p\wedge \lfa_q\wedge \lfa_r$
has 2-termed representation in the form of a Clifford algebraic polynomial. The same number of
terms is expected for each of (\ref{note:check})
in the form of a quaternionic polynomial. To obtain such a form we need to not only embed
$\qh$ into $CL({\mathbb R}^4)$, but also project $CL({\mathbb R}^4)$ onto $\qh$; they are the two
maps $\iota$ and $\pi$ below.

\bn
Denote by $H$ the 4-dimensional linear subspace of $CL(\mathbb{R}^4)$ spanned by
$1,\e_{23},\e_{13},\e_{12}$. Define a linear isomorphism $\iota: H\mapsto {\mathbb H}$:
\be
\iota(1)=1, \ \ \
\iota(\e_{23})=-\i, \ \ \
\iota(\e_{13})=\j, \ \ \
\iota(\e_{12})=-\k.
\ee
Let $\pi_H$ be the restriction of $CL(\mathbb{R}^4)$ to subspace $H$. Define
\be\ba{rll}
\pi: \hspace{.2cm} CL(\mathbb{R}^4) &\rightarrow& \mathbb{H} \\
                        \A &\mapsto   & (\iota\circ \pi_H) (\A(1+\e_0)(1-\I_4))
\ea
\label{def:pi}
\ee
\en

\vskip .2cm
While $H=\langle 1, \e_{23},\e_{13},\e_{12} \rangle_{\mathbb R}$, we have
\be\ba{lll}
\e_0 H &=&  \langle \e_0, \e_{023},\e_{013},\e_{012} \rangle_{\mathbb R}, \\

\I_4 H &=&  \langle \e_{0123}, \e_{01},\e_{02},\e_{03} \rangle_{\mathbb R}, \\

\e_0 \I_4 H &=& \langle \e_{123}, \e_{1},\e_{2},\e_{3} \rangle_{\mathbb R}.
\ea\ee
Then $CL(\mathbb{R}^4)=H\oplus \e_0 H\oplus \I_4 H\oplus \e_0 \I_4 H$.
By definition, 
\be\ba{lrrrll}
  \{   1, &     \e_0,  &  \e_{123},   &  -\e_{0123}\}  &\xrightarrow[]{\pi}& 1;  \\
  \{\e_1, & -\e_{01},  &  -\e_{23},   &  -\e_{023}\}   &\xrightarrow[]{\pi}& \i; \\
  \{\e_2, & -\e_{02},  &   \e_{13},   &   \e_{013}\}   &\xrightarrow[]{\pi}& \j; \\
  \{\e_3, & -\e_{03},  &  -\e_{12},   &  -\e_{012}\}   &\xrightarrow[]{\pi}& \k.
\ea
\label{ker:pi}
\ee

For any $\a\in \qh$, obviously $\pi(\lfa)=\a$, so $\pi\circ \f={\rm id}$ in $\qh$.
Furthermore, 
$\pi$ is a linear isomorphism from any of $H, \e_0 H, \I_4 H, \e_0 \I_4 H$ to $\mathbb{H}$.

Let $K_-$ be the odd elements of $ker(\pi)\subset CL(\mathbb{R}^4)$. 
By (\ref{ker:pi}), $K_-$ is a 4-dimensional real space spanned by 
\[\ba{lll}
\e_0-\e_{123} &=& \e_0(1-\I_4), \\

\e_1+\e_{023} &=& \e_1(1-\I_4), \\

\e_2-\e_{013} &=& \e_2(1-\I_4), \\

\e_3+\e_{012} &=& \e_3(1-\I_4).
\ea
\]
So 
\be
K_-={\mathbb R}^4(1-\I_4)=\langle \e_l(1-\I_4),\ l=0..3\rangle_{\mathbb R}.
\label{ker:odd}
\ee

%{\it Remark}: Before defining $\pi$, if we require $\pi$ to be a 
%linear isomorphism from any of $H, \e_0 H, \I_4 H, \e_0 \I_4 H$ to $\mathbb{H}$, such that
%$\pi(1)=-\pi(\I_4)=1$ and $\pi(\e_1)=\i$, then by letting 
%$\pi(\A)=(\iota\circ \pi_H) (\A(1+\lambda\e_0)(1+\mu\I_4))$ for any 
%$\A\in CL(\mathbb{R}^4)$, we get $\lambda=-\mu=1$, which is just (\ref{def:pi}).

\bp \label{conjugate homomorphism:prop}
Let $\A,\B\in CL(\mathbb{R}^4)$, then
\be \label{conjugate homomorphism}
\pi(\A\B)=\left\{
  \begin{array}{ll} \vspace{.1cm}
    \pi(\A)\pi(\B),                         & \hbox{if~$\A$~is~even;} \\
    \vspace{.1cm}
    \pi(\A)\pi(\overline{\B}),   & \hbox{if~$\A$~is~odd.}
  \end{array}
\right.
\ee
In particular, $\pi(\lfa\lfb) = \a\overline{\b}$ for any $\a, \b\in \qh$.
\ep

\bo
(1). When $\A$ is even, let $\A=A_0+\A_v+D_0\I_4+\D_v\I_4$,
where $A_0,D_0$ are scalars, and $\A_v,\D_v\in\langle \e_{23},\e_{13},\e_{12} \rangle_{\mathbb R}$.
Then $\pi(\A)=\iota(A_0+\A_v-D_0-\D_v)$.

If $\B$ is even, let $\B=B_0+\B_v+C_0\I_4+\C_v\I_4$, where $B_0,C_0$ are scalars,
and $\B_v,\C_v\in\langle \e_{23},\e_{13},\e_{12} \rangle_{\mathbb R}$. Then
$\pi(\B)=\iota(B_0+\B_v-C_0-\C_v)$. By
\begin{eqnarray*}
% \nonumber to remove numbering (before each equation)
  \A\B &=& A_0B_0+A_0\B_v+A_0C_0\I_4+A_0\C_v\I_4 \\
       & & +B_0\A_v+\A_v\B_v+C_0\A_v\I_4+\A_v\C_v\I_4 \\
       & & +B_0D_0\I_4+D_0\I_4\B_v+C_0D_0+D_0\C_v \\
       & & +B_0\D_v\I_4+\D_v\B_v\I_4+C_0\D_v+\D_v\C_v,
\end{eqnarray*}
we have
\begin{eqnarray*}
% \nonumber to remove numbering (before each equation)
  \pi(\A\B) &=& \iota(A_0B_0+A_0\B_v-A_0C_0-A_0\C_v\\
         & &   +B_0\A_v+\A_v\B_v-C_0\A_v-\A_v\C_v \\
       & & +C_0D_0+D_0\C_v-B_0D_0-D_0\B_v \\
       & & -B_0\D_v-\D_v\B_v+C_0\D_v+\D_v\C_v) \\
       &=& \pi(\A)\iota(B_0+\B_v-C_0-\C_v) \\
       &=& \pi(\A)\pi(\B).
\end{eqnarray*}

If $\B$ is odd, let
$\B=b_0\e_0+\b_v+c_0\e_0\I_4+\c_v\I_4$, where $b_0,c_0$ are scalars,
and $\b_v,\c_v\in\langle \e_1,\e_2,\e_3 \rangle_{\mathbb R}$. Then
$\pi(\B)=\iota(b_0-\b_v\e_0\I_4+c_0-\c_v\e_0\I_4)$. By
\begin{eqnarray*}
% \nonumber to remove numbering (before each equation)
  \A\B &=& A_0b_0\e_0+A_0\b_v+A_0c_0\e_0\I_4+A_0\c_v\I_4 \\
       & & +b_0\A_v\e_0+\A_v\b_v+c_0\A_v\e_0\I_4+\A_v\c_v\I_4 \\
       & & -b_0D_0\e_0\I_4+D_0\I_4\b_v-c_0D_0\e_0-D_0\c_v \\
       & & -b_0\D_v\e_0\I_4+\D_v\I_4\b_v-c_0\D_v\e_0-\D_v\c_v,
    \end{eqnarray*}
we have
\begin{eqnarray*}
% \nonumber to remove numbering (before each equation)
  \pi(\A\B) &=& \iota(A_0b_0-A_0\c_v\e_0\I_4-A_0\b_v\e_0\I_4+A_0c_0\\
       & & +b_0\A_v-\A_v\c_v\e_0\I_4-\A_v\b_v\e_0\I_4+c_0\A_v \\
       & & +D_0\b_v\e_0\I_4-c_0D_0-b_0D_0+D_0\c_v\e_0\I_4 \\
       & & +\D_v\b_v\e_0\I_4-c_0\D_v-b_0\D_v+\D_v\c_v\e_0\I_4) \\
       &=& \pi(\A)\iota(b_0-\b_v\e_0\I_4+c_0-\c_v\e_0\I_4) \\
       &=& \pi(\A)\pi(\B).
\end{eqnarray*}

(2). When $\A$ is odd, let $\A=a_0\e_0+\a_v+d_0\e_0\I_4+\d_v\I_4$, where $a_0,d_0$ are
scalars, and $\a_v,\d_v\in\langle \e_1,\e_2,\e_3 \rangle_{\mathbb R}$. Then
$\pi(\A)=\iota(a_0-\a_v\e_0\I_4+d_0-\d_v\e_0\I_4)$.

If $\B$ is even, let $\B=B_0+\B_v+C_0\I_4+\C_v\I_4$, where $B_0,C_0$ are scalars,
and $\B_v,\C_v\in\langle \e_{23},\e_{13},\e_{12} \rangle_{\mathbb R}$. Then
$\overline{\B}=B_0+\B_v-C_0\I_4-\C_v\I_4$, and $\pi(\overline{\B})=\iota(B_0+\B_v+C_0+\C_v)$.
By
\begin{eqnarray*}
% \nonumber to remove numbering (before each equation)
  \A\B &=& a_0B_0\e_0+a_0\e_0\B_v+a_0C_0\e_0\I_4+a_0\C_v\e_0\I_4 \\
       & & +B_0\a_v+\a_v\B_v+C_0\a_v\I_4+\a_v\C_v\I_4 \\
       & & +d_0 B_0\e_0\I_4+d_0\e_0\I_4\B_v+d_0C_0\e_0+d_0\C_v\e_0 \\
       & & +B_0\d_v\I_4+\d_v\I_4\B_v+C_0\d_v+\d_v\C_v,
\end{eqnarray*}
we have
\begin{eqnarray*}
% \nonumber to remove numbering (before each equation)
  \pi(\A\B) &=& \iota(a_0B_0+a_0\B_v+a_0C_0+a_0\C_v \\
       & & -B_0\a_v\e_0\I_4-\a_v\B_v\e_0\I_4-C_0\a_v\e_0\I_4-\a_v\C_v\e_0\I_4 \\
       & & +d_0 B_0+d_0\B_v+d_0C_0+d_0\C_v \\
       & & -B_0\d_v\e_0\I_4-\d_v\B_v\e_0\I_4-C_0\d_v\e_0\I_4-\d_v\C_v\e_0\I_4) \\
       &=& \pi(\A)\iota(B_0+\B_v+C_0+\C_v) \\
       &=& \pi(\A)\pi(\overline{\B}).
\end{eqnarray*}

If $\B$ is odd, let $\B=b_0\e_0+\b_v+\c_v\I_4+c_0\e_0\I_4$, where $b_0,c_0$ are scalars and
$\b_v,\c_v\in\langle \e_1,\e_2,\e_3 \rangle_{\mathbb R}$. Then
$\overline{\B}=b_0\e_0-\b_v+\c_v\I_4-c_0\e_0\I_4$, and
$\pi(\overline{\B})=\iota(b_0+\b_v\e_0\I_4-c_0-\c_v\e_0\I_4)$.
By
\begin{eqnarray*}
% \nonumber to remove numbering (before each equation)
  \A\B &=& a_0b_0+a_0\e_0\b_v+a_0c_0\I_4-a_0\c_v\e_0\I_4 \\
       & & +b_0\a_v\e_0+\a_v\b_v+c_0\a_v\e_0\I_4+\a_v\c_v\I_4 \\
       & & -d_0b_0\I_4+d_0\e_0\I_4\b_v-d_0c_0+d_0\c_v\e_0 \\
       & & -b_0\d_v\e_0\I_4+\d_v\I_4\b_v-c_0\d_v\e_0-\d_v\c_v,
\end{eqnarray*}
we have
\begin{eqnarray*}
% \nonumber to remove numbering (before each equation)
  \pi(\A\B) &=& a_0b_0-a_0\c_v\e_0\I_4+a_0\b_v\e_0\I_4-a_0c_0\\
         & & -b_0\a_v\e_0\I_4+\a_v\b_v+c_0\a_v\e_0\I_4 -\a_v\c_v \\
       & &  +d_0b_0-d_0\c_v\e_0\I_4+d_0\e_0\I_4\b_v-d_0c_0 \\
       & & -b_0\d_v\e_0\I_4-\d_v\c_v+\d_v\b_v+c_0\d_v\e_0\I_4 \\
       &=& \pi(\A)\iota(b_0+\b_v\e_0\I_4-c_0-\c_v\e_0\I_4) \\
       &=& \pi(\A)\pi(\overline{\B}).
\end{eqnarray*}
\eo

\bc
For any $r>0$ and $\a_l\in \qh$,
\be \label{conjugate homomorphism:vector}
\pi(\lfa_1\lfa_2\cdots \lfa_r)
=\left\{\ba{lll}
\a_1\bar{\a}_2\a_3\bar{\a}_4\cdots \a_{r-1}\bar{\a}_r, & \hbox{if} & r \hbox{ is even}; \\
\a_1\bar{\a}_2\a_3\bar{\a}_4\cdots \bar{\a}_{r-1}\a_r, & \hbox{if} & r \hbox{ is odd}.
\ea\right.
\ee
\ec

\bp
For any $\a_l\in \qh$,
\be\label{det:expr}
[\lfa_1\lfa_2\lfa_3\lfa_4]=-\frac{1}{4}
(\a_1\bar{\a}_2\a_3\bar{\a}_4
+\a_4\bar{\a}_3\a_2\bar{\a}_1
-\a_4\bar{\a}_1\a_2\bar{\a}_3
-\a_3\bar{\a}_2\a_1\bar{\a}_4).
\ee
\ep

\bo
By definition,
$
\lfa_1\wedge\lfa_2\wedge\lfa_3\wedge\lfa_4=[\lfa_1\lfa_2\lfa_3\lfa_4]\I_4.
$
On the other hand, by Lemma \ref{lem:ungrade},
$4\lfa_1\wedge\lfa_2\wedge\lfa_3\wedge\lfa_4 
=\lfa_1\lfa_2\lfa_3\lfa_4+\lfa_4\lfa_3\lfa_2\lfa_1
-\lfa_4\lfa_1\lfa_2\lfa_3-\lfa_3\lfa_2\lfa_1\lfa_4$.
So
\begin{eqnarray*}
% \nonumber to remove numbering (before each equation)
{[}\lfa_1\lfa_2\lfa_3\lfa_4] &=& 
\pi(\lfa_1\wedge\lfa_2\wedge\lfa_3\wedge\lfa_4)\, \pi(\I_4) \\
&=& -\frac{1}{4}(\a_1\bar{\a}_2\a_3\bar{\a}_4
+\a_4\bar{\a}_3\a_2\bar{\a}_1-\a_4\bar{\a}_1\a_2\bar{\a}_3
-\a_3\bar{\a}_2\a_1\bar{\a}_4).
\end{eqnarray*}
\eo

\bp For any $\a_1, \a_2, \a_3\in \qh$,
\be\label{expr:dual:1}
\a_{123}^\sim=\frac{1}{2}(\a_1\bar{\a}_2\a_3-\a_3\bar{\a}_2\a_1).
\ee
Furthermore,
\be
\bar{\a}_{123}^{\sim}=-\overline{\a_{123}^{\sim}}.
\label{expr:dual:2}
\ee
\ep

\bo By Lemma \ref{lem:ungrade},
$2(\lfa_1\wedge\lfa_2\wedge\lfa_3)^\sim =
(\lfa_1\lfa_2\lfa_3-\lfa_3\lfa_2\lfa_1)\I_4$.
So
\[
2 \a_{123}^\sim = 2 \pi(\lfa_1\wedge\lfa_2\wedge\lfa_3)\, \pi(\overline{\I_4}) 
               =  \pi(\lfa_1\lfa_2\lfa_3-\lfa_3\lfa_2\lfa_1) 
               = \a_1\bar{\a}_2\a_3-\a_3\bar{\a}_2\a_1.
\]
Similarly, $2 \bar{\a}_{123}^\sim=\bar{\a}_1\a_2\bar{\a}_3
-\bar{\a}_3\a_2\bar{\a}_1
=-2\overline{\a_{123}^{\sim}}$.
\eo

\section{The associated real linear system approach}
\setcounter{equation}{0}

The revised starting form (\ref{rev:start}) can be written as a real linear system 
$\A\lfq=\lfd$, where matrix $\A$ is given by
(\ref{express:A}). The solution is $\lfq=adj(\A)\lfd/\det(\A)$ under the non-degeneracy condition 
$\det(\A)\neq 0$.

\bt
\be
\det(\A) = 2\hskip -.12cm \sum_{i,j=0}^3
(\lfa_i\cdot\lfa_j)^2
-\Big(\sum_{i=0}^3
\lfa_i\cdot\lfa_i\Big)^2\hskip -.1cm
-8[\lfa_0\lfa_1\lfa_2\lfa_3].
\label{expr:deta}
\ee
\et

\bo
After direct expansion and simplification by \verb"Maple", we get that $\det(\A)$ is a 
polynomial of 196 terms in the indeterminates $a_{ij}$ for
$i,j=0..3$. The terms can be categorized into two groups. The first group $G_1$
is the sum of terms of the form $a_{ij}^4$ or
$a_{uv}^2a_{kl}^2$; the second group $G_2$ is the sum of the other terms.
All terms of $G_2$ are square-free.

We first consider $G_1$.
Denote
\[
\N_1=\left(
              \begin{array}{cccc}
                1 & 1 & 1 & 1 \\
                1 & 1 & 1 & 1 \\
                1 & 1 & 1 & 1 \\
                1 & 1 & 1 & 1 \\
              \end{array}
            \right),
\hskip .6cm
\N_2=\left(
              \begin{array}{cccc}
                \hfill 1 & -1       & -1       & -1       \\
                -1       & \hfill 1 & -1       & -1       \\
                -1       & -1       & \hfill 1 & -1       \\
                -1       & -1       & -1       & \hfill 1 \\
              \end{array}
\right).
\]
Further denote
\[
\v_l^T:=(a_{l0}^2,a_{l1}^2,a_{l2}^2,a_{l3}^2),\ \hbox{ for }l=0..3.
\]
It is easy to prove that for any $i,j\in \{0,1,2,3\}$,
\[
\v_i^T\N_1\v_j=\lfa_i^2 \lfa_j^2, \hskip .4cm
\v_i^T\N_2\v_j=2\v_i^T\v_j-\v_i^T\N_1\v_j.
\]

By direct verification, we confirm that $G_1$ equals
\be\ba{ll}
& (\v_0^T\ \v_1^T\ \v_2^T\ \v_3^T)\left(
               \begin{array}{cccc}
                  \N_1 & \N_2 & \N_2 & \N_2 \\
                  \N_2 & \N_1 & \N_2 & \N_2 \\
                  \N_2 & \N_2 & \N_1 & \N_2 \\
                  \N_2 & \N_2 & \N_2 & \N_1 \\
               \end{array}
             \right)
             \left(\ba{l}
             \v_0 \\
             \v_1\\
             \v_2\\
             \v_3
             \ea\right)
             \\
             [2mm]
=& \ds \sum_{i=0}^3\v_i^T\N_1\v_i+2\hskip -.2cm \sum_{0\leq j<k\leq 3}\hskip -.2cm  \v_j^T\N_2\v_k
\\
[2mm]
=& \ds \sum_{i=0}^3 \lfa_i^2
-2\hskip -.2cm \sum_{0\leq j<k\leq 3}\hskip -.2cm  \lfa_j^2 \lfa_k^2
+4\hskip -.2cm \sum_{0\leq j<k\leq 3}\hskip -.2cm  \v_j^T\v_k.
\ea
\label{group:1}
\ee

Again
by direct verification, we confirm that $G_2$ equals
\be
-8[\lfa_0 \lfa_1 \lfa_2 \lfa_3]+4\sum_{0\leq i<j\leq 3}\hskip -.2cm
(\lfa_i\cdot\lfa_j)^2
-4\sum_{0\leq i<j\leq 3}\hskip -.2cm
\v_i^T\v_j.
\label{h:expr}
\ee  
By (\ref{group:1}) and (\ref{h:expr}), we get (\ref{expr:deta}).
\eo

Next consider the expression of the adjugate $adj(\A)$ of $\A$.
Denote
\be\ba{rllll}
\Lambda &:=& \ds \sum_{l=0}^3 \lfa_l^2 &=& \ds \sum_{l=0}^3\a_l\overlinea_l, \\

adj(\a_i) &:=& \ds
\sum_{l=0}^3 (\lfa_i\cdot \lfa_l) \overlinea_l
&=&\ds \frac{1}{2}\sum_{l=0}^3 (\a_i\overlinea_l+\a_l\overlinea_i) \overlinea_l.\Bigstrut
\ea
\label{note:lambda}
\ee

\bt
With the notations introduced by (\ref{note:check}) and (\ref{note:lambda}),
\be
\ba{lll}
adj(\A)
&=&
\phantom{-}\big(
\phantom{-}\ \, 2 \overlinea_{123}^{\sim}+2 adj(\overlinea_0)-\Lambda\, 
\overlinea_0\left|\right.1
\big)\\
[.5mm]
&&
-\big(
-2 \overlinea_{023}^{\sim}+2 adj(\overlinea_1)-\Lambda\, \overlinea_1\left|\right. \i\,
\big)\\
[.5mm]
&&
-\big(\phantom{-}\ \,
2 \overlinea_{013}^{\sim}+2 adj(\overlinea_2)-\Lambda\, \overlinea_2\left|\right. \j\,
\big)\\
[.5mm]
&&
-\big(
-2 \overlinea_{012}^{\sim}+2 adj(\overlinea_3)-\Lambda\, \overlinea_3\left|\right. \k
\big).
\ea
\label{expr:adja}
\ee
\et

\bo The entries of $adj(\A)$ are algebraic minors of $\A$. Once we obtain them, then
by Proposition \ref{matrix-quaternion-form},
\be
adj(\A)=(\p_0\,|\,1)+(\p_1\,|\,\i)+(\p_2\,|\,\j)+(\p_3\,|\,\k),
\label{adj:express}
\ee
where the $\p_l$ satisfy (\ref{soln:a4}) for $(m_{ij})_{i,j=0..3}=adj(\A)$.

We first compute $\p_0=p_{00}+p_{01}\i+p_{02}\j+p_{03}\k$. By (\ref{soln:a4}) and direct verification,
\[\hskip -.12cm
\ba{lll}
% \nonumber to remove numbering (before each equation)
  p_{00} 
         &=&a_{00}\left(\lfa_0^2-\lfa_1^2-\lfa_2^2-\lfa_3^2\right)
         +2a_{10}\lfa_0\cdot\lfa_1+2a_{20}\lfa_0\cdot\lfa_2\\
         
         && \hfill +2a_{30}\lfa_0\cdot\lfa_3
          -2[\e_0\lfa_1\lfa_2\lfa_3],
\\ [1mm]

  p_{01} &
          =& -a_{01}
          \left(\lfa_0^2-\lfa_1^2-\lfa_2^2-\lfa_3^2\right)
         -2a_{11}\lfa_0\cdot\lfa_1-2a_{21}\lfa_0\cdot\lfa_2\\
         
         && \hfill -2a_{31}\lfa_0\cdot\lfa_3
          +2[\e_1\lfa_1\lfa_2\lfa_3],
\\ [1mm]

  p_{02} &
         =&
         \left(\lfa_0^2-\lfa_1^2-\lfa_2^2-\lfa_3^2\right)
         -2a_{12}\lfa_0\cdot\lfa_1-2a_{22}\lfa_0\cdot\lfa_2\\
         
         && \hfill -2a_{32}\lfa_0\cdot\lfa_3
          +2[\e_2 \lfa_1 \lfa_2 \lfa_3],
\\ [1mm]

  p_{03} &
         =& -a_{03}
         \left(\lfa_0^2-\lfa_1^2-\lfa_2^2-\lfa_3^2\right)
         -2a_{13}\lfa_0\cdot\lfa_1-2a_{23}\lfa_0\cdot\lfa_2\\
         
         && \hfill -2a_{33}\lfa_0\cdot\lfa_3
          +2[\e_3 \lfa_1 \lfa_2 \lfa_3].
\ea\]
So
\be\ba{rl}
  \p_0 =&\hskip -.2cm \left(\lfa_0^2-\lfa_1^2-\lfa_2^2-\lfa_3^2\right)
  \lfba_0
  +2\left(\lfa_0\cdot\lfa_1\right)\lfba_1
  +2\left(\lfa_0\cdot\lfa_2\right)\lfba_2
  +2\left(\lfa_0\cdot\lfa_3\right)\lfba_3 \\
[1mm]
&\hskip -.2cm -2\big([\e_0 \lfa_1 \lfa_2 \lfa_3]\e_0
-[\e_1 \lfa_1 \lfa_2 \lfa_3]\e_1
-[\e_2 \lfa_1 \lfa_2 \lfa_3]\e_2
-[\e_3 \lfa_1 \lfa_2 \lfa_3]\e_3\big).
\ea\ee

Similarly, we get
\be\ba{rl}
  \p_1 =&\hskip -.2cm \left(\lfa_0^2-\lfa_1^2+\lfa_2^2+\lfa_3^2\right)\lfba_1
  -2\left(\lfa_1\cdot\lfa_0\right)\lfba_0
  -2\left(\lfa_1\cdot\lfa_2\right)\lfba_2
  -2\left(\lfa_1\cdot\lfa_3\right)\lfba_3 \\
[1mm]
&\hskip -.2cm -2\big([\e_0 \lfa_0 \lfa_2 \lfa_3]\e_0
-[\e_1 \lfa_0 \lfa_2 \lfa_3]\e_1
-[\e_2 \lfa_0 \lfa_2 \lfa_3]\e_2
-[\e_3 \lfa_0 \lfa_2 \lfa_3]\e_3\big),\\
[2.5mm]

  \p_2 =&\hskip -.2cm \left(\lfa_0^2+\lfa_1^2-\lfa_2^2+\lfa_3^2\right)\lfba_2
  +2\left(\lfa_2\cdot\lfa_0\right)\lfba_0
  +2\left(\lfa_2\cdot\lfa_1\right)\lfba_1
  +2\left(\lfa_2\cdot\lfa_3\right)\lfba_3 \\
[1mm]
&\hskip -.2cm +2\big([\e_0 \lfa_0 \lfa_1 \lfa_3]\e_0
-[\e_1 \lfa_0 \lfa_1 \lfa_3]\e_1
-[\e_2 \lfa_0 \lfa_1 \lfa_3]\e_2
-[\e_3 \lfa_0 \lfa_1 \lfa_3]\e_3\big),\\
[2.5mm]

  \p_3 =&\hskip -.2cm \left(\lfa_0^2-\lfa_1^2-\lfa_2^2-\lfa_3^2\right)\lfba_3
  +2\left(\lfa_3\cdot\lfa_0\right)\lfba_0
  +2\left(\lfa_3\cdot\lfa_2\right)\lfba_1
  +2\left(\lfa_3\cdot\lfa_3\right)\lfba_2 \\
[1mm]
&\hskip -.2cm -2\big([\e_0 \lfa_0 \lfa_1 \lfa_2]\e_0
-[\e_1 \lfa_0 \lfa_1 \lfa_2]\e_1
-[\e_2 \lfa_0 \lfa_1 \lfa_2]\e_2
-[\e_3 \lfa_0 \lfa_1 \lfa_2]\e_3\big).
\ea\ee
Substituting these expressions into (\ref{adj:express}), and using the
notations introduced by (\ref{note:check}) and  (\ref{note:lambda}), we get
(\ref{expr:adja}).
\eo

Now consider the general linear quaternionic equation (\ref{eqn:general}). Let
$\b_p=b_{p0}+b_{p1}\i+b_{p2}\j+b_{p3}\k$ for $p=1..n-1$, and set
\be
\a_i:=\sum_{p=1}^{n-1} b_{pi} \c_p, \hbox{ for } i\in \{0,1,2,3\}.
\label{a:realform}
\ee
Then (\ref{eqn:general}) is changed into the revised starting form (\ref{rev:start}).

\bt \label{general-solution:thm}
With the $\a_i$ taking values (\ref{a:realform}), 
\be\ba{rl}
\det(\A) =&\hskip -.2cm \ds 2\hskip -.2cm
\sum_{p,q,r,s=1}^{n-1}\hskip -.2cm
(\lfc_p\cdot\lfc_q)(\lfc_r\cdot\lfc_s)
(\lfb_p\cdot\lfb_r)(\lfb_q\cdot\lfb_s) \\
%[2mm]
&\hskip -.2cm \ds\hfill -\Big( \sum_{p,q=1}^{n-1} 
(\lfc_p\cdot\lfc_q)(\lfb_p\cdot\lfb_q)\Big)^2
-\frac{1}{3} \sum_{p,q,r,s=1}^{n-1}
[\lfc_p \lfc_q \lfc_r \lfc_s]\,[\lfb_p \lfb_q \lfb_r \lfb_s].
\ea
\label{deta:last}
\ee
The solution of (\ref{eqn:general}) is unique if and only if $\det(\A)\neq 0$, and
the solution when $\det(\A)\neq 0$ is
\be
\q=\frac{-1}{3\det(\A)}\sum_{p,q,r=1}^{n-1}\hskip -.15cm
\Big\{\overline{\c}_{pqr}^{\sim} \d \overline{\b}_{pqr}^{\sim}
+3(\lfc_p\cdot\lfc_q)\Big((\lfb_p\cdot\lfb_q) \overline{\c}_r\d\overline{\b}_r
-2(\lfb_p\cdot\lfb_r)\overline{\c}_r\d\overline{\b}_q \Big)\Big\}.
\label{general-solution}
\ee
\et

\bo We have $\lfa_i\cdot\lfa_j=\sum_{p,q=1}^{n-1}(\c_p\cdot\c_q)b_{pi}b_{qj}$. So
\be\ba{lll}
\ds \sum_{i,j=0}^3(\lfa_i\cdot\lfa_j)^2
   &\hskip -.2cm =& \ds\hskip -.3cm \sum_{i,j=0}^3\ \sum_{p,q,r,s=1}^{n-1}  \hskip -.1cm
   b_{pi}b_{qj}b_{ri}b_{sj} (\lfc_p\cdot\lfc_q)(\lfc_r\cdot\lfc_s) \\
   [1mm]
   &\hskip -.2cm =&\hskip -.3cm \ds \sum_{p,q,r,s=1}^{n-1}\hskip -.1cm
   (\lfc_p\cdot\lfc_q)(\lfc_r\cdot\lfc_s)(\lfb_p\cdot\lfb_r)(\lfb_q\cdot\lfb_s);
\\ [2mm]

\ds \Big(\sum_{i=0}^3 \lfa_i\cdot \lfa_i\Big)^2
&\hskip -.2cm =&\hskip -.2cm  \ds \Bigg(\sum_{i=0}^3\
\sum_{p,q,r=1}^{n-1} b_{pi}b_{qi}\lfc_p\cdot\lfc_q\Bigg)^2 \\
&\hskip -.2cm =&\hskip -.2cm  \ds \Bigg(\sum_{p,q=1}^{n-1}
(\lfc_p\cdot\lfc_q)(\lfb_p\cdot\lfb_q)\Bigg)^2;
\\
[2mm]
\ds {[}\lfa_1 \lfa_2 \lfa_3 \lfa_4]
&\hskip -.2cm =&\hskip -.2cm  \ds  \sum_{p,q,r,s=1}^{n-1}
b_{p1}b_{q2}b_{r3}b_{s4}[\lfc_p \lfc_q \lfc_r \lfc_s] \\
[2mm]
&\hskip -.2cm =&\hskip -.4cm  \ds \sum_{1\leq p<q<r<s\leq n-1}
\hskip -.4cm
[\lfc_p \lfc_q \lfc_r \lfc_s]
[\lfb_p \lfb_q \lfb_r \lfb_s].
\ea
\ee
Substituting them into (\ref{expr:deta}), we get (\ref{deta:last}).

Consider the action of $adj(\A)$ upon $\d$. We divide the 
expression of $adj(\A)$ in (\ref{expr:adja}) into two parts,
the first part $P_1$ is $
(2 \overlinea_{123}^{\sim}\,|\,1)
+(2 \overlinea_{023}^{\sim}\,|\,\i)
-(2 \overlinea_{013}^{\sim}\,|\,\j)
+(2 \overlinea_{012}^{\sim}\,|\,\k)$, and the second part $P_2$ is the rest.

For $P_1$, 
\be\ba{rl}
\overlinea_{123}^{\sim} & \hskip -.2cm
=\ds
\sum_{i=0}^3\ [\e_i \lfba_1 \lfba_2 \lfba_3]\f^{-1}(\e_i)\\

& \hskip -.2cm =\ds
\sum_{i=0}^3\
\sum_{p,q,r=1}^{n-1} [\e_i \lfbc_p \lfbc_q \lfbc_r] b_{p1}b_{q2}b_{r3} \f^{-1}(\e_i)\\

& \hskip -.2cm =\ds
\sum_{i=0}^3\
\sum_{1\leq p<q<r\leq n-1}
[\e_i \lfbc_p \lfbc_q \lfbc_r][\e_0 \lfb_p \lfb_q \lfb_r]\f^{-1}(\e_i).
\ea\ee
Similarly, if denoting 
\[
\overlinea_{\check{1}}^{\sim}:=\overlinea_{023}^{\sim},\ \ \
\overlinea_{\check{2}}^{\sim}:=\overlinea_{013}^{\sim},\ \ \
\overlinea_{\check{3}}^{\sim}:=\overlinea_{012}^{\sim},
\]
then for $j=1,2,3$,
\be
\overlinea_{\check{j}}^{\sim} =
\sum_{i=0}^3\
\sum_{1\leq p<q<r\leq n-1} (-1)^{j-1}
[\e_i \lfbc_p  \lfbc_q  \lfbc_r]
[\e_j \lfb_p \lfb_q \lfb_r]\f^{-1}(\e_i).
\ee
So
\[\ba{rl}
& \overlinea_{123}^{\sim}\d
+\overlinea_{023}^{\sim}\d\,\i
-\overlinea_{013}^{\sim}\d\,\j
+\overlinea_{012}^{\sim}\d\,\k
\\[2mm]

=& \ds
\sum_{i,j=0}^3\
\sum_{1\leq p<q<r\leq n-1}
[\e_i \lfbc_p  \lfbc_q  \lfbc_r]
[\e_j \lfb_p \lfb_q \lfb_r]\f^{-1}(\e_i)\,\d\, \overline{\f^{-1}(\e_j)} \\[2mm]

=& \ds
-\sum_{i,j=0}^3\
\sum_{1\leq p<q<r\leq n-1}
[\e_i \lfbc_p  \lfbc_q  \lfbc_r]
[\e_j \lfbb_p \lfbb_q  \lfbb_r]\f^{-1}(\e_i)\,\d\, \f^{-1}(\e_j) \\[5mm]

=& \ds -\sum_{1\leq p<q<r\leq n-1} \overline{\c}_{pqr}^{\sim} \d\, \overline{\b}_{pqr}^{\sim}
\\[5mm]

=& \ds
-\frac{1}{6}\sum_{p,q,r=1}^{n-1} \overline{\c}_{pqr}^{\sim} \d\, \overline{\b}_{pqr}^{\sim}.
\ea
\]

For $P_2$, we have
\[\ba{rl}
& (2 adj(\overlinea_0)-\Lambda\, \overlinea_0)\d
-(2 adj(\overlinea_1)-\Lambda\, \overlinea_1)\d\i\\

& \hfill
-(2 adj(\overlinea_2)-\Lambda\, \overlinea_2)\d\j
-(2 adj(\overlinea_3)-\Lambda\, \overlinea_3)\d\k
\\
[2mm]

=& \ds
\sum_{p,q,r=1}^{n-1}(\lfc_p\cdot\lfc_q)\overline{\c}_r\d
\Big( (b_{p0}b_{q0}-b_{p1}b_{q1}-b_{p2}b_{q2}-b_{p3}b_{q3})b_{r0}\\

& \hskip 2.65cm +(b_{p0}b_{q0}-b_{p1}b_{q1}+b_{p2}b_{q2}+b_{p3}b_{q3})b_{r1}\i\\ [1.5mm]

& \hskip 2.65cm +(b_{p0}b_{q0}+b_{p1}b_{q1}-b_{p2}b_{q2}+b_{p3}b_{q3})b_{r2}\j\\

& \hskip 2.65cm +(b_{p0}b_{q0}+b_{p1}b_{q1}+b_{p2}b_{q2}-b_{p3}b_{q3})b_{r3}\k\Big)\\
[3mm]

& \ds +2\sum_{p,q,r=1}^{n-1}(\lfc_p\cdot\lfc_q)\overline{\c}_r\d\Big(
b_{p0}b_{q0}b_{r0}-b_{p0}b_{r0}(b_{q1}\i+b_{q2}\j+b_{q3}\k)\\

& \hskip 3.15cm +b_{p0}b_{q1}b_{r1}-b_{p1}b_{q1}b_{r1}\i-b_{p1}b_{q2}b_{r1}\j-b_{p1}b_{q3}b_{r1}\k\\
[1.5mm]

& \hskip 3.15cm +b_{p0}b_{q2}b_{r2}-b_{p1}b_{q2}b_{r2}\i-b_{p2}b_{q2}b_{r2}\j-b_{p2}b_{q3}b_{r2}\k\\

& \hskip 3.15cm +b_{p0}b_{q3}b_{r3}-b_{p1}b_{q3}b_{r3}\i-b_{p2}b_{q3}b_{r3}\j-b_{p3}b_{q3}b_{r3}\k
\Big) \\
[3mm]

=& \ds -\sum_{p,q,r=1}^{n-1}(\c_p\cdot\c_q)(\b_p\cdot\b_q)\overline{\c}_r\d\overline{\b}_r
+2\ds\sum_{p,q,r=1}^{n-1}(\c_p\cdot\c_q)(\b_p\cdot\b_r)\overline{\c}_r\d\overline{\b}_q.
\ea\]

Combining the results of $P_1, P_2$, we get (\ref{general-solution}).
\eo

{\it Remark}: the solution (\ref{general-solution}) is just the following identity
in $\qh$:
\be\ba{lll}
& \ds \sum_{p,q,r=1}^{n-1}\Big\{(\overline{\c}_{pqr}^{\sim}
\mid\overline{\b}_{pqr}^{\sim})+3(\lfc_p\cdot\lfc_q)
\Big((\lfb_p\cdot\lfb_q)(\overline{\c}_r\mid\overline{\b}_r) \\
& \hfill\ds -2(\lfb_p\cdot\lfb_r)(\overline{\c}_r\mid\overline{\b}_q) \Big) \Big\}
\Big(\sum_{s=1}^{n-1} (\c_s\,|\,\b_s)\Big) \\
&= \ds -3\det\big(\sum_{p=1}^{n-1} (\c_p\,|\,\b_p)\big)(1\,|\,1).
\ea
\ee

\section{The Clifford algebraic approach}

Assume that the input linear quaternionic equation
(\ref{eqn:general})
is the image of an equation in $CL({\mathbb R}^4)$ under the projection
$\pi$, where the preimages of the $\c_i, \q, \b_i, \d$ under $\pi$ are 
vectors of ${\mathbb R}^4$. 
The corresponding equation in $CL({\mathbb R}^4)$ is called the {\it lift} of the 
quaternionic equation. 

\bl For any $\b_i, \c_i, \d, \q\in \qh$,
$\sum_{i=1}^{n-1}\c_i\q\b_i-\d=0$ if and only if
\be
(\sum_{i=1}^{n-1}\lfc_i\overline{\lfq}\lfb_i-\lfd)(1+\I_4)=0.
\label{linear Clifford equation:1}
\ee
\el 

\bo
Obviously 
$\pi((\sum_{i=1}^{n-1}\lfc_i\overline{\lfq}\lfb_i-\lfd)(1+\I_4))
=(\sum_{i=1}^{n-1}\c_i\q\b_i-\d)(1+\pi(\overline{\I_4}))
=2(\sum_{i=1}^{n-1}\c_i\q\b_i-\d)$. So 
(\ref{linear Clifford equation:1}) leads to the input equation
(\ref{eqn:general}).

Conversely, by (\ref{ker:odd}), for an $\A\in K_-$ (the odd elements of $ker(\pi)$),
there exists an $\x\in {\mathbb R}^4$ such that $\A=\x(1-\I_4)$. 
Since $\pi(\sum_{i=1}^{n-1}\lfc_i\overline{\lfq}\lfb_i-\lfd)=
\sum_{i=1}^{n-1}\c_i\q\b_i-\d=0$, 
there exists $\x\in {\mathbb R}^4$ such that
\be
\sum_{i=1}^{n-1}\lfc_i\overline{\lfq}\lfb_i-\lfd = \x(1-\I_4).
\ee
Multiplying both sides from the right by $1+\I_4$, and using
$(1-\I_4)(1+\I_4)=1-\I_4^2=0$, we get  
(\ref{linear Clifford equation:1}).
\eo

Below we solve (\ref{linear Clifford equation:1})
for vector variable
$\overline{\lfq}$.

Let $I_{init}$ be the two-sided ideal generated by  
$\sum_{i=1}^{n-1}\c_i\q\b_i-\d\in \qh$.
Let $I_{+}$ be the $\mathbb R$-linear combination of the following elements
in $CL({\mathbb R}^4)$:
\be
\B (\sum_{i=1}^{n-1}\lfc_i\overline{\lfq}\lfb_i-\lfd)(1+\I_4) \C, 
\ee 
where $\B, \C\in CL({\mathbb R}^4)$ and $\B$ is an even element.
$I_{+}$ is called the {\it left-even two-sided ideal}
generated by $(\sum_{i=1}^{n-1}\lfc_i\overline{\lfq}\lfb_i-\lfd)(1+\I_4)$.

\bl
$\pi(I_{+})\subset I_{init}$.
\el

\bo The result follows
\[
\pi\Big(\B (\sum_{i=1}^{n-1}\lfc_i\overline{\lfq}\lfb_i-\lfd)(1+\I_4) \C\Big)
=2 \pi(\B) (\sum_{i=1}^{n-1}\c_i\q\b_i-\d) \pi(\overline{\C})
\in I_{init}.
\]
\eo

So
if $\A\in CL({\mathbb R}^4)$ is an odd element, then
$\e_0\A (\sum_{i=1}^{n-1}\lfc_i\overline{\lfq}\lfb_i-\lfd)(1+\I_4) \C
\in I_{+}$, and
$\pi\Big(\e_0\A (\sum_{i=1}^{n-1}\lfc_i\overline{\lfq}\lfb_i-\lfd)(1+\I_4) \C\Big)
\in I_{init}$.

\textbf{Notation.} For $\b_i, \c_j\in \qh$, denote
\be\ba{ll}
\lfb_{ij}:=\lfb_{i}\wedge\lfb_{j}, & \lfb_{ijk}:=\lfb_{i}\wedge\lfb_{j}\wedge\lfb_{k},\\
\lfc_{ij}:=\lfc_{i}\wedge\lfc_{j}, & \lfc_{ijk}:=\lfc_{i}\wedge\lfc_{j}\wedge\lfc_{k}.
\ea
\ee

Applying the linear operator 
$\sum_{p,q,r=1}^{n-1} (\e_0\lfc_{pqr}^\sim \mid \lfb_{pqr}^\sim)$ 
to both sides of (\ref{linear Clifford equation:1}), we get
\be \label{linear Clifford equation:11}
\sum_{i,p,q,r=1}^{n-1}
\e_0\lfc_{pqr}^\sim\lfc_i\overline{\lfq}\lfb_i\lfb_{pqr}^\sim(1-\I_4)
=\sum_{p,q,r=1}^{n-1}\e_0\lfc_{pqr}^\sim\lfd\lfb_{pqr}^\sim(1-\I_4).
\ee

\bl \cite{hestenes}
\label{lem:hestenes}
For any $\x_l\in {\mathbb R}^4$,
\be\ba{rll}
\x_1\x_2 &=& \phantom{-} \x_1\cdot \x_2+\x_1\wedge \x_2, \\[1mm]

\x_1(\x_2\wedge \x_3\wedge \x_4)^\sim &=& \phantom{-} 
[\x_1\x_2\x_3\x_4]+(\x_1\cdot \x_2)(\x_3\wedge \x_4)^\sim \\
&& \hfill
-(\x_1\cdot \x_3)(\x_2\wedge \x_4)^\sim
+(\x_1\cdot \x_4)(\x_1\wedge \x_4)^\sim, \\[1mm]

(\x_2\wedge \x_3\wedge \x_4)^\sim \x_1 &=& \phantom{-} 
[\x_1\x_2\x_3\x_4]-(\x_1\cdot \x_2)(\x_3\wedge \x_4)^\sim\\
&& \hfill 
+(\x_1\cdot \x_3)(\x_2\wedge \x_4)^\sim
-(\x_1\cdot \x_4)(\x_1\wedge \x_4)^\sim.
\ea
\ee
\el

Applying Lemma \ref{lem:hestenes} to
 the left side of (\ref{linear Clifford equation:11}), we get
\begin{eqnarray*}
% \nonumber to remove numbering (before each equation)
  & & L.H.S. \\
  &=& \sum_{i,p,q,r=1}^{n-1} \e_0
  \Big([\lfc_i\lfc_p\lfc_q\lfc_r]-((\lfc_i\cdot\lfc_p)\lfc_{qr}
  -(\lfc_i\cdot\lfc_q)\lfc_{pr}+(\lfc_i\cdot\lfc_r)\lfc_{pq})^\sim\Big)\overline{\lfq} \\
  & &  \Big([\lfb_i\lfb_p\lfb_q\lfb_r]+((\lfb_i\cdot\lfb_p)\lfb_{qr}
  -(\lfb_i\cdot\lfb_q)\lfc_{pr}+(\lfb_i\cdot\lfb_r)\lfb_{pq})^\sim\Big)(1-\I_4) \\
  &:=& \Sigma_1+\Sigma_2-\Sigma_3-\Sigma_4,\Bigstrut
\end{eqnarray*}
where
\be\ba{lll}
\Sigma_1 &\hskip -.2cm =&\hskip -.5cm  \ds\sum_{i,p,q,r=1}^{n-1}[\lfc_i\lfc_p\lfc_q\lfc_r][\lfb_i\lfb_p\lfb_q\lfb_r]
\e_0\overline{\lfq}(1-\I_4); 
\\
\Sigma_2 &\hskip -.2cm =&\hskip -.5cm \ds\sum_{i,p,q,r=1}^{n-1}[\lfc_i\lfc_p\lfc_q\lfc_r]\e_0\overline{\lfq}
((\lfb_i\cdot\lfb_p)\lfb_{qr}
-(\lfb_i\cdot\lfb_q)\lfb_{pr}+(\lfb_i\cdot\lfb_r)\lfb_{pq})(\I_4-1); \Bigstrut
\\
\Sigma_3 &\hskip -.2cm =&\hskip -.5cm \ds\sum_{i,p,q,r=1}^{n-1}[\lfb_i\lfb_p\lfb_q\lfb_r]\e_0((\lfc_i\cdot\lfc_p)\lfc_{qr}
-(\lfc_i\cdot\lfc_q)\lfc_{pr}+(\lfc_i\cdot\lfc_r)\lfc_{pq})
\overline{\lfq}(1-\I_4);\Bigstrut
\\
\Sigma_4 &\hskip -.2cm =&\hskip -.5cm \ds\sum_{i,p,q,r=1}^{n-1} \e_0((\lfc_i\cdot\lfc_p)\lfc_{qr}-(\lfc_i\cdot\lfc_q)\lfc_{pr}
+(\lfc_i\cdot\lfc_r)\lfc_{pq})\overline{\lfq} \Bigstrut\\
&&\hfill
         ((\lfb_i\cdot\lfb_p)\lfb_{qr}
         -(\lfb_i\cdot\lfb_q)\lfb_{pr}+(\lfb_i\cdot\lfb_r)\lfb_{pq})(\I_4-1).
\ea\ee

By the symmetry of the inner product and the antisymmetry of the outer product,
$\Sigma_2=\Sigma_3=0$. For $\Sigma_4$, again by the two symmetries, 
\be
\Sigma_4 = \sum_{i,p,q,r=1}^{n-1}(\lfc_i\cdot\lfc_p)
\Big(
-3(\lfb_i\cdot\lfb_p)\e_0\lfc_{qr}\overline{\lfq}
\lfb_{qr}
+6(\lfb_i\cdot\lfb_q)\e_0\lfc_{qr}
  \overline{\lfq}\lfb_{pr}\Big)(1-\I_4).
\ee
Thus we get the following result:

\bl
The following element is in $I_{+}$: 
\be\ba{l}
 \ds \left(
-\sum_{p,q,r=1}^{n-1}\e_0\lfc_{pqr}^\sim\lfd\lfb_{pqr}^\sim
+
\sum_{i,p,q,r=1}^{n-1}[\lfc_i\lfc_p\lfc_q\lfc_r][\lfb_i\lfb_p\lfb_q\lfb_r]\e_0\overline{\lfq} 
\right.\\
 \ds 
 \left.
+\sum_{i,p,q,r=1}^{n-1}(\lfc_i\cdot\lfc_p)
\Big(
6(\lfb_i\cdot\lfb_q)\e_0\lfc_{qr}
  \overline{\lfq}\lfb_{pr}
  -3(\lfb_i\cdot\lfb_p)\e_0\lfc_{qr}\overline{\lfq}
\lfb_{qr}\Big)
\right)(1-\I_4).
\ea
\label{lem:step1}
\ee
\el

In (\ref{lem:step1}), the last summation can be generated by applying two
linear operators to the left side of equation 
(\ref{linear Clifford equation:1}) respectively: 
applying linear operator $\sum_{i,p,q=1}^{n-1}(\lfc_i\cdot\lfc_p)
(\lfb_i\cdot\lfb_q)(\e_0\lfc_q \mid \lfb_p)$ to
(\ref{linear Clifford equation:1}), we get
\be
\label{linear Clifford equation:12} 
\ba{l}
\ds
\sum_{i,p,q,r=1}^{n-1}(\lfc_i\cdot\lfc_p)(\lfb_i\cdot\lfb_q)\e_0\lfc_q\lfc_r
\overline{\lfq}\lfb_r\lfb_p(1-\I_4) \\

\hskip 2cm \ds
=\sum_{i,p,q=1}^{n-1}(\lfc_i\cdot\lfc_p)(\lfb_i\cdot\lfb_q)\e_0\lfc_q\lfd\lfb_p(1-\I_4);
\ea
\ee
applying $\sum_{i,p,q=1}^{n-1}(\lfc_i\cdot\lfc_p)(\lfb_i\cdot\lfb_p)
(\e_0\lfc_q \mid \lfb_q)$ to (\ref{linear Clifford equation:1}), we get
\be
\label{linear Clifford equation:13}
\ba{l}
\ds
\sum_{i,p,q,r=1}^{n-1}(\lfc_i\cdot\lfc_p)(\lfb_i\cdot\lfb_p)
\e_0\lfc_q\lfc_r\overline{\lfq}\lfb_r\lfb_q(1-\I_4)\\

\hskip 2cm \ds
=\sum_{i,p,q=1}^{n-1}(\lfc_i\cdot\lfc_p)(\lfb_i\cdot\lfb_p)\lfc_q\lfd\lfb_q(1-\I_4).
\ea
\ee

\bl 
The following two elements are in $I_{+}$: 
\be \label{linear Clifford equation:12:left}
\ba{l} 
\ds  \left(
-\sum_{i,p,q=1}^{n-1}(\lfc_i\cdot\lfc_p)(\lfb_i\cdot\lfb_q)\e_0\lfc_q\lfd\lfb_p \right.\\
\ds
+\sum_{i,p,q,r=1}^{n-1}(\lfc_i\cdot\lfc_p)(\lfb_i\cdot\lfb_q)
(\lfc_q\cdot\lfc_r)(\lfb_p\cdot\lfb_r)\e_0\overline{\lfq}\\

\hfill \ds \left.-\sum_{i,p,q,r=1}^{n-1}(\lfc_i\cdot\lfc_p)(\lfb_i\cdot\lfb_q)
\e_0\lfc_{qr}\lfq\lfb_{pr}\right)(1-\I_4);
\ea\ee
\be\ba{l} \label{linear Clifford equation:13:left}
\ds  
\left(-\sum_{i,p,q=1}^{n-1}(\lfc_i\cdot\lfc_p)(\lfb_i\cdot\lfb_p)\e_0\lfc_q\lfd\lfb_q
\right.\\
\ds
+\sum_{i,p,q,r=1}^{n-1}(\lfc_i\cdot\lfc_p)(\lfb_i\cdot\lfb_p)
(\lfc_q\cdot\lfc_r)(\lfb_q\cdot\lfb_r)\e_0\overline{\lfq} \\
\hfill \ds \left.-\sum_{i,p,q,r=1}^{n-1}(\lfc_i\cdot\lfc_p)(\lfb_i\cdot\lfb_p)
\e_0\lfc_{qr}\lfq\lfb_{qr}\right)(1-\I_4).
\ea\ee
\el

\bo
Applying Lemma \ref{lem:hestenes} to
 the left side of (\ref{linear Clifford equation:12}), then by the 
 symmetry of the inner product and the antisymmetry of the outer product,
 we get
\begin{eqnarray*}
% \nonumber to remove numbering (before each equation)
L.H.S.
   &=& \sum_{i,p,q,r=1}^{n-1}(\lfc_i\cdot\lfc_p)(\lfb_i\cdot\lfb_q)(\lfc_q\cdot\lfc_r
   +\lfc_{qr})\e_0\overline{\lfq} (\lfb_r\cdot\lfb_p-\lfb_{pr})(1-\I_4) \\
   &=& \sum_{i,p,q,r=1}^{n-1}(\lfc_i\cdot\lfc_p)(\lfb_i\cdot\lfb_q)
   (\lfc_q\cdot\lfc_r)(\lfb_r\cdot\lfb_p)\e_0\overline{\lfq}(1-\I_4) \\
   & & -\sum_{i,p,q,r=1}^{n-1}(\lfc_i\cdot\lfc_p)(\lfb_i\cdot\lfb_q)(\lfc_q\cdot\lfc_r)
   \e_0\overline{\lfq}\lfb_{pr}(1-\I_4) \\
   & & +\sum_{i,p,q,r=1}^{n-1}(\lfc_i\cdot\lfc_p)(\lfb_i\cdot\lfb_q)(\lfb_r\cdot\lfb_p)
   \e_0\lfc_{qr}\overline{\lfq}(1-\I_4) \\
   & & -\sum_{i,p,q,r=1}^{n-1}(\lfc_i\cdot\lfc_p)(\lfb_i\cdot\lfb_q)
   \e_0\lfc_{qr}\overline{\lfq}\lfb_{pr}(1-\I_4)
   \\
   &=& \sum_{i,p,q,r=1}^{n-1}(\lfc_i\cdot\lfc_p)(\lfb_i\cdot\lfb_q)
   (\lfc_q\cdot\lfc_r)(\lfb_r\cdot\lfb_p)\e_0\overline{\lfq}(1-\I_4) \\
   & & -\sum_{i,p,q,r=1}^{n-1}(\lfc_i\cdot\lfc_p)(\lfb_i\cdot\lfb_q)
   \e_0\lfc_{qr}\overline{\lfq}\lfb_{pr}(1-\I_4).
\end{eqnarray*} 
From the left side of (\ref{linear Clifford equation:13}), we get
a similar expression. They lead to the results of the lemma.
\eo

By the above lemmas, we get

\bp The following element is in $I_{+}$: 
\be\ba{l} \label{general-solution:clifford}
\ds \sum_{p,q,r,s=1}^{n-1}\Big\{[\lfc_p\lfc_q\lfc_r\lfc_s][\lfb_p\lfb_q\lfb_r\lfb_s] 
+3(\lfc_p\cdot\lfc_q)(\lfc_r\cdot\lfc_s)(\lfb_p\cdot\lfb_q)(\lfb_r\cdot\lfb_s) \\
\hfill \ds -6(\lfc_p\cdot\lfc_q)(\lfc_r\cdot\lfc_s)
(\lfb_q\cdot\lfb_s)(\lfb_p\cdot\lfb_r)\Big\}\e_0\overline{\lfq}(1-\I_4) \\
\ds -\sum_{p,q,r=1}^{n-1}\Big\{\e_0\lfc_{pqr}^\sim\lfd\lfb_{pqr}^\sim
+3(\lfc_p\cdot\lfc_r)(\lfb_p\cdot\lfb_r)\e_0\lfc_q\lfd\lfb_q \\
\hfill -6(\lfc_p\cdot\lfc_r)(\lfb_q\cdot\lfb_r)\e_0\lfc_q\lfd \lfb_p\Big\}(1-\I_4).
\ea\ee
\ep

\bc \label{lem:nc}
The following element is in $I_{init}$:

\be\ba{l} \label{general-solution:quaternion}
\ds \sum_{p,q,r,s=1}^{n-1}\Big\{[\lfc_p\lfc_q\lfc_r\lfc_s][\lfb_p\lfb_q\lfb_r\lfb_s]
+3(\lfc_p\cdot\lfc_q)(\lfc_r\cdot\lfc_s)(\lfb_p\cdot\lfb_q)(\lfb_r\cdot\lfb_s) \\
\hfill 
-6(\lfc_p\cdot\lfc_q)(\lfc_r\cdot\lfc_s)(\lfb_q\cdot\lfb_s)(\lfb_p\cdot\lfb_r
)\Big\}\q \\
\ds -\sum_{p,q,r=1}^{n-1}\Big\{\overline{\c}_{pqr}^\sim\d\overline{\b}_{pqr}^\sim
+3(\lfc_p\cdot\lfc_r)(\lfb_p\cdot\lfb_r)\overline{\c}_q\d\overline{\b}_q
-6(\lfc_p\cdot\lfc_r)(\lfb_q\cdot\lfb_r)\overline{\c}_q\d\overline{\b}_p\Big\}.
\ea\ee
It leads to Theorem \ref{general-solution:thm} directly.
\ec

\vskip .2cm
Set
\be\ba{lll}
\Delta &\hskip -.2cm :=&\ds\hskip -.2cm 
\sum_{p,q,r,s=1}^{n-1}\Big\{[\lfc_p\lfc_q\lfc_r\lfc_s][\lfb_p\lfb_q\lfb_r\lfb_s]
+3(\lfc_p\cdot\lfc_q)(\lfc_r\cdot\lfc_s)(\lfb_p\cdot\lfb_q)(\lfb_r\cdot\lfb_s) \\

&& \ds\hfill 
-6(\lfc_p\cdot\lfc_q)(\lfc_r\cdot\lfc_s)(\lfb_q\cdot\lfb_s)(\lfb_p\cdot\lfb_r
)\Big\}, \\

\Phi &\hskip -.2cm :=&\hskip -.2cm  \ds \sum_{p,q,r=1}^{n-1}
\Big\{(\overline{\c}_{pqr}^{\sim}\mid\overline{\b}_{pqr}^{\sim})
+3(\lfc_p\cdot\lfc_r)
(\lfb_p\cdot\lfb_r)(\overline{\c}_q\mid\overline{\b}_q)\\

&& \ds\hfill 
-6(\lfc_p\cdot\lfc_r)
(\lfb_q\cdot\lfb_r)(\overline{\c}_q\mid\overline{\b}_p)
\Big\}.
\ea
\label{NOT:PHI}
\ee
Then $\Delta\in \mathbb R$, and
$\Phi$ is a linear operator in $\qh$.
Let $\A\q=\d$ be the associated real linear system of the
input quaternionic equation, then by (\ref{deta:last}),
$\Delta=-3\det(\A)$.
By Corollary \ref{lem:nc}, 
$\Phi\A\q=\Phi\d=\Delta\, \q=-3\det(\A) \q$. So in matrix form,
\be
\Phi=-3\, adj(\A).
\ee
The non-degeneracy of $\Phi$ is the same with that of $\A$.

\section{General linear quaternionic equation with conjugate}

A general linear quaternionic equation containing the conjugate of the
quaternionic indeterminate $\q$ is of the following form:
\be \label{linear H equation:conjugate case}
\sum_{r=1}^\alpha \c_r\q\b_r-\sum_{r=\alpha+1}^{n-1} \c_r\overline{\q}\b_r = \d,
\ee
where $0\leq \alpha<n$, and the $\c_r, \b_r, \d$ are given quaternions.

Let $\q=x_0+\x$, where $x_0, \x$ are the real part and the pure imaginary 
part of $\q$ respectively. Set 
\be
\h:=\sum_{r=1}^{\alpha} \c_r\b_r-\sum_{r=\alpha+1}^{n-1} \c_r\b_r. 
\ee
Then (\ref{linear H equation:conjugate case}) can be written as
\be \label{conj:eqn}
\sum_{r=1}^{n-1}\c_r\x\b_r=\d-x_0\h.
\ee

Applying Theorem \ref{general-solution:thm} to equation (\ref{conj:eqn}), 
then using the notation of (\ref{NOT:PHI}), we get
\be
\Delta\,\x = \Phi\d-x_0\Phi\h.
\label{conj:1}
\ee 
By $\textmd{Re}(\x)=0$, we have $\textmd{Re}(\Phi\d)-x_0\textmd{Re}(\Phi\h)=0$, 
so 
\be
x_0=\frac{\textmd{Re}(\Phi\d)}{\textmd{Re}(\Phi\h)}.
\label{conj:2}
\ee

\bl
Let $\M$ be 
the coefficient matrix of the associated real linear system of 
(\ref{linear H equation:conjugate case}), then
\be
\det(\M)=-\frac{1}{3}\textmd{Re}(\Phi\h).
\label{last:M}
\ee
\el

\bo
By definition, 
\[
\M\left(
1,
\i,
\j,
\k
\right)
=\left(
\h,
\A\i,
\A\j,
\A\k
\right),
\]
where $\A$ is the coefficient matrix of the associated real linear equations of 
(\ref{eqn:general}). When $\A$ is invertible,
\[
\A^{-1}\M\left(
1,
\i,
\j,
\k
\right)
=\left(
\A^{-1}\h,
\i,
\j,
\k
\right)
=\left(
\textmd{Re}(\A^{-1}\h) + \textmd{Im} (\A^{-1}\h),
\i,
\j,
\k
\right),
\]
so 
\[\ba{lll}
\det(\A^{-1}\M) &=& (\det(\A))^{-1}\det(\M)\\ [1mm]

 &=& \textmd{Re}(\A^{-1}\h)\\[1mm]
 
 &=& (\det(\A))^{-1} \textmd{Re}(adj(\A)\h)\\[1mm]
 
  &=& \ds-\frac{1}{3} (\det(\A))^{-1} \textmd{Re}(\Phi\h).
\ea\]

When $\A$ is not invertible, (\ref{last:M}) is still valid by the continuity
of both sides of it in the coordinates of the quaternionic coefficients
of (\ref{linear H equation:conjugate case}).
\eo

The following result is direct from (\ref{conj:1}), (\ref{conj:2}) and (\ref{last:M}):  

\bt
When $\det(\A)=-\Delta/3\neq 0$ and $\det(\M)=-\textmd{Re}(\Phi\h)/3\neq 0$,
the solution to (\ref{linear H equation:conjugate case}) is
\be
\Delta\,\textmd{Re}(\Phi\h)\,\q
= \Delta\,\textmd{Re}(\Phi\d)-\textmd{Re}(\Phi\d)\, \Phi\h
+\textmd{Re}(\Phi\h)\,\Phi\d.
\ee
\et

%\section{Conclusion}

%In this paper we solve the long last problem of finding explicit basis-free 
%solutions of general linear quaternionic equations
%by the Grassmann algebra approach. The technique might be extended to linear 
%quaternionic equations in one quaternionic variable
%and its conjugate, or more than one quaternionic variables.

\end{document}